\documentclass{article}

\usepackage[dcucite]{harvard}

\usepackage{amssymb}
\usepackage{amsmath}

\newtheorem{thm}{Theorem}
\newtheorem{defn}[thm]{Definition}
\newtheorem{rem}[thm]{Remark}
\newenvironment{pf}[1][Proof]{\textbf{#1.} }{\ \rule{0.5em}{0.5em} \bigskip}

\newenvironment{keywords}{\begin{center}
\begin{minipage}[c]{11cm} {\bf Keywords:}} {\end{minipage}
\end{center}}
\newenvironment{msc}{\begin{center}
\begin{minipage}[c]{11cm} {\bf 2000 Mathematics Subject Classification:}} {\end{minipage}
\end{center} \bigskip}

\begin{document}

\title{Noether's Theorem for Fractional\\ Optimal Control Problems\footnote{To be
presented at FDA'06 -- \emph{2nd IFAC Workshop on Fractional
Differentiation and its Applications}, 19-21 July 2006, Porto,
Portugal. Accepted (07-March-2006) for the Conference Proceedings.
Research Report CM06/I-11.}}

\author{Gast\~{a}o S. F. Frederico\thanks{Supported by IPAD
(Instituto Portugu\^{e}s de Apoio ao Desenvolvimento).}\\
\texttt{gfrederico@mat.ua.pt}
\and
Delfim F. M. Torres\thanks{Supported by CEOC (Centre for Research on Optimization
and Control) through FCT (Portuguese Foundation for Science and Technology),
cofinanced by the European Community fund FEDER/POCTI.}\\
\texttt{delfim@mat.ua.pt}}

\date{Department of Mathematics\\
University of Aveiro\\
3810-193 Aveiro, Portugal}

\maketitle

\begin{abstract}
We begin by reporting on some recent results of the authors \cite{tncdf},
concerning the use of the fractional Euler-Lagrange notion to
prove a Noether-like theorem for the problems of the calculus of variations
with fractional derivatives. We then obtain, following the Lagrange
multiplier technique used in \cite{Agrawal:2004a}, a new version
of Noether's theorem to fractional optimal control systems.
\end{abstract}


\begin{msc}
49K05, 26A33, 70H33.
\end{msc}


\section{Introduction}

The concept of symmetry plays an important role both in Physics and
Mathematics. Symmetries are described by transformations of the
system, which result in the same object after the transformation
is carried out. They are described mathematically by parameter
groups of transformations. Their importance ranges from fundamental
and theoretical aspects to concrete applications, having profound
implications in the dynamical behavior of the systems, and in
their basic qualitative properties.

Another fundamental notion in Physics and Mathematics is the one
of conservation law. Typical application of conservation laws in
the calculus of variations and optimal control is to reduce the
number of degrees of freedom, and thus reducing the problems to a
lower dimension, facilitating the integration of the differential
equations given by the necessary optimality conditions.

Emmy Noether was the first who proved, in 1918, that the notions of
symmetry and conservation law are connected: when a system
exhibits a symmetry, then a conservation law can be obtained. One
of the most important and well known illustration of this deep
and rich relation, is given by the conservation of energy in
Mechanics: the autonomous Lagrangian $L(q,\dot{q})$, correspondent
to a mechanical system of conservative points, is invariant under
time-translations (time-homogeneity symmetry), and\footnote{We use
the notation $\partial_i f$ to denote the partial derivative
of some function $f$ with respect to its $i$-th argument.}
\begin{equation}
\label{eq:R}
\begin{gathered}
-L(q,\dot{q})+\partial_{2} L(q,\dot{q})\cdot\dot{q}
\equiv \text{constant}
\end{gathered}
\end{equation}
follows from Noether's theorem,
\textrm{i.e.}, the total energy of a conservative closed system
always remain constant in time, ``it cannot be created or
destroyed, but only transferred from one form into another''.
Expression \eqref{eq:R} is valid along all the
Euler-Lagrange extremals $q(\cdot)$ of an autonomous problem of
the calculus of variations. The conservation law
\eqref{eq:R} is known in the calculus of variations as the
2nd Erdmann necessary condition; in concrete applications, it gains
different interpretations: conservation of energy in Mechanics;
income-wealth law in Economics; first law of Thermodynamics; etc.
The literature on Noether's theorem is vast, and many extensions
of the classical results of Emmy Noether are now available for the
more general setting of optimal control (see \cite{delfimEJC,delfimPortMath04}
and references therein). Here we remark that in all those results
conservation laws always refer to problems with integer derivatives.

Nowadays fractional differentiation plays an important role in
various fields: physics (classic and quantum mechanics,
thermodynamics, etc), chemistry, biology, economics, engineering,
signal and image processing, and control theory
\cite{Agrawal:2004b,CD:Hilfer:2000,CD:Klimek:2002}. Its origin goes back more
than 300 years, when in 1695 L'Hopital asked Leibniz the
meaning of $\frac{d^{n}y}{dx^{n}}$ for $n=\frac{1}{2}$. After
that, many famous mathematicians, like J.~Fourier,
N.~H.~Abel, J.~Liouville, B.~Riemann, among others, contributed to
the development of the Fractional Calculus
\cite{CD:Hilfer:2000,CD:MilRos:1993,CD:SaKiMa:1993}.

The study of fractional problems of
the Calculus of Variations and
respective Euler-Lagrange type equations is a subject
of current strong research.
F.~Riewe \cite{CD:Riewe:1996,CD:Riewe:1997} obtained a version of
the Euler-Lagrange equations for problems of the Calculus of
Variations with fractional derivatives, that combines the
conservative and non-conservative cases. In 2002 O.~Agrawal
proved a formulation for variational problems with right and left
fractional derivatives in the Riemann-Liouville sense
\cite{CD:Agrawal:2002}. Then these Euler-Lagrange equations were
used by D.~Baleanu and T.~Avkar to investigate problems with
Lagrangians which are linear on the velocities
\cite{CD:BalAv:2004}. In \cite{Klimek2001,MR1966935} fractional problems of the calculus of
variations with symmetric fractional derivatives are considered and correspondent
Euler-Lagrange equations obtained, using both Lagrangian and Hamiltonian
formalisms. In all the above mentioned studies,
Euler-Lagrange equations depend on left and right
fractional derivatives, even when the problem depend only on
one type of them. In \cite{Klimek2005} problems depending on
symmetric derivatives are considered for which
Euler-Lagrange equations include only the
derivatives that appear in the formulation of the
problem. In \cite{El-Nabulsi2005b,El-Nabulsi2005a} Riemann-Liouville
fractional integral functionals, depending on a parameter $\alpha$ but not on
fractional-order derivatives of order $\alpha$, are introduced and respective
fractional Euler-Lagrange type equations obtained. More recently, the authors have used
the results of \cite{CD:Agrawal:2002} to generalize Noether's theorem for the
context of the Fractional Calculus of Variations \cite{tncdf}.
Here we extend the previous optimal control Noether results
in \cite{delfimEJC,delfimPortMath04}
to the wider context of fractional optimal control, making use
(i) of the fractional version of Noether's theorem obtained
by the authors in \cite{tncdf}, (ii) and the Lagrange multiplier
rule \cite{Agrawal:2004a}.


\section{Fractional derivatives}
\label{sec:fdRL}

In this section we collect the definitions of right and left
Riemann-Liouville fractional derivatives and their main properties
\cite{CD:Agrawal:2002,CD:MilRos:1993,CD:SaKiMa:1993}.

\begin{defn}
Let $f$ be a continuous and integrable function in the interval
$[a,b]$. For all $t \in [a,b]$, the left Riemann-Liouville
fractional derivative $_aD_t^\alpha f(t)$, and the right
Riemann-Liouville fractional derivative $_tD_b^\alpha f(t)$, of
order $\alpha$, are defined in the following way:
\begin{equation}
\label{eq:DFRLE}
_aD_t^\alpha f(t) =
\frac{1}{\Gamma(n-\alpha)}\left(\frac{d}{dt}\right)^{n}
\int_a^t (t-\theta)^{n-\alpha-1}f(\theta)d\theta \, ,
\end{equation}
\begin{equation}
\label{eq:DFRLD}
_tD_b^\alpha f(t) =
\frac{1}{\Gamma(n-\alpha)}\left(-\frac{d}{dt}\right)^{n}
\int_t^b(t-\theta)^{n-\alpha-1}f(\theta)d\theta \, ,
\end{equation}
where $n \in \mathbb{N}$, $n-1 \leq \alpha < n$, and $\Gamma$ is
the Euler gamma function.
\end{defn}

\begin{rem}
If $\alpha$ is an integer, then from \eqref{eq:DFRLE}
and \eqref{eq:DFRLD} one obtains the standard derivatives, that is,
\begin{gather*}
\label{eq:DU}
_aD_t^\alpha f(t) = \left(\frac{d}{dt}\right)^\alpha f(t) \, , \\
_tD_b^\alpha f(t) = \left(-\frac{d}{dt}\right)^\alpha f(t) \, .
\end{gather*}
\end{rem}

\begin{thm} Let $f$ and $g$ be two continuous functions
on $[a,b]$. Then, for all $t \in [a,b]$, the following properties
hold:
\begin{enumerate}
\item for $p>0$, $$_aD_t^p\left(f(t)+g(t)\right)
= {_aD_t^p}f(t)+{_aD_t^p}g(t) \, ;$$
\item for $p \geq q \geq 0$,
$$_aD_t^p\left(_aD_t^{-q} f(t)\right) = {_aD_t^{p-q}}f(t) \, ;$$
\item for $p>0$, $$_aD_t^p\left(_aD_t^{-p} f(t)\right) = f(t)$$
(fundamental property of the Riemann-Liouville fractional derivatives);
\item for $p>0$, $$\int_a^b \left(_aD_t^p f(t)\right)g(t)dt
= \int_a^b f(t) \, _tD_b^p g(t)dt \, .$$
\end{enumerate}
\end{thm}

\begin{rem}
In general, the fractional derivative of a constant is not
equal to zero.
\end{rem}

\begin{rem}
The fractional derivative of order $p>0$
of function $(t-a)^\upsilon$, $\upsilon>-1$,
is given by
\begin{equation*}
_aD_t^p(t-a)^\upsilon
= \frac{\Gamma(\upsilon+1)}{\Gamma(-p+\upsilon+1)}(t-a)^{\upsilon-p} \, .
\end{equation*}
\end{rem}

\begin{rem}
In the literature, when one reads ``Riemann-Liouville fractional derivative'',
one usually means the ``left Riemann-Liouville fractional derivative''.
In Physics, if $t$ denotes the time-variable, the
right Riemann-Liouville fractional derivative of $f(t)$
is interpreted as a future state of the process $f(t)$. For this
reason, the right-derivative is usually
neglected in applications, when the present state of the process does
not depend on the results of the future development. Following \cite{Agrawal:2004a},
in this work we only consider problems with left Riemann-Liouville fractional derivatives.
Using  \cite{tncdf}, the results of the paper can, however, be written for the case
when both left and right fractional derivatives are present.
\end{rem}

We refer the reader interested in additional background
on fractional theory, to the comprehensive book
\cite{CD:SaKiMa:1993}.


\section{Preliminaries}
\label{sec:MR}

In \cite{CD:Agrawal:2002} a formulation of the
Euler-Lagrange equations is given for problems of the calculus of
variations with fractional derivatives.

Let us consider the following fractional problem
of the calculus of variations: to find function
$q(\cdot)$ that minimizes the integral functional
\begin{gather}
\label{Pf}
I[q(\cdot)] = \int_a^b L\left(t,q(t),{_aD_t^\alpha} q(t)\right) dt \, ,
\end{gather}
where the Lagrangian
$L :[a,b] \times \mathbb{R}^{n} \times
\mathbb{R}^{n} \rightarrow \mathbb{R}$
is a $C^{2}$ function with respect to all its arguments, and
$0 < \alpha \leq 1$.

\begin{rem}
In the case $\alpha = 1$, problem \eqref{Pf} is reduced to
the classical problem
\begin{equation*}
I[q(\cdot)] = \int_a^b L\left(t,q(t),\dot{q}(t)\right) dt
\longrightarrow \min \, .
\end{equation*}
\end{rem}

\begin{thm}[\cite{CD:Agrawal:2002}]
\label{Thm:FractELeq} If $q$ is a minimizer of problem \eqref{Pf},
then it satisfies the \emph{fractional Euler-Lagrange equations}:
\begin{equation}
\label{eq:eldf}
\partial_{2} L\left(t,q,{_aD_t^\alpha q}\right) \\
+ {_tD_b^\alpha}\partial_{3} L\left(t,q,{_aD_t^\alpha q}\right) = 0 \, .
\end{equation}
\end{thm}

The following definition is useful in order to introduce an
appropriate concept of \emph{fractional conservation law}.

\begin{defn}[\cite{tncdf}]
Given two functions $f$ and $g$ of class $C^{1}$ in the interval
$[a,b]$, we introduce the following notation:
\begin{equation*}
\mathcal{D}\left\{(f g)_{t}^{\alpha}\right\}
= -g \, {_tD_b^\alpha} f + f \, {_aD_t^\alpha} g \, ,
\end{equation*}
where $t \in [a,b]$.
\end{defn}

\begin{rem}
For $\alpha = 1$ operator $\mathcal{D}$ is reduced to
\begin{equation*}
\begin{split}
\mathcal{D}\left\{(f g)_{t}^{1}\right\}
&=- g \, {_tD_b^1 f} + f \, {_aD_t^1} g \\
&= \dot{f} g + f \dot{g} = \frac{d}{dt}(f g) \, .
\end{split}
\end{equation*}
\end{rem}

\begin{rem}
The linearity of the operators $_aD_t^\alpha$ and $_tD_b^\alpha$
imply the linearity of the operator $\mathcal{D}$.
\end{rem}

\begin{defn}[\cite{tncdf}]
\label{eq:fcl}
We say that $C_{f}\left(t,q,{_aD_t^\alpha q}\right)$, where
$C_{f}$ has the form of a sum of products
\begin{equation}
\label{eq:somaPrd}
C_{f}\left(t,q,d\right)
= \sum_{i} C_{i}^1\left(t,q,d\right) \cdot C_{i}^2\left(t,q,d\right)
\end{equation}
is a \emph{fractional conservation law} if, and only if,
\begin{equation}
\label{eq:def:lcf}
\mathcal{D}\left\{{C}_{f}\left(t,q,{_aD_t^\alpha q}\right)\right\} = 0
\end{equation}
along all the fractional Euler-Lagrange extremals
(\textrm{i.e.} along all the solutions of the fractional
Euler-Lagrange equations \eqref{eq:eldf}).
\end{defn}

\begin{rem}
For $\alpha = 1$ \eqref{eq:def:lcf} is reduced to
\begin{equation*}
\frac{d}{dt}\left\{{C}_{f}\left(t,q(t),\dot{q}(t)\right)\right\} = 0 \\
\Leftrightarrow {C}_{f}\left(t,q(t),\dot{q}(t)\right) \equiv \text{constant} \, ,
\end{equation*}
which is the standard meaning of \emph{conservation law}:
a function $C_{f}\left(t,q,\dot{q}\right)$ preserved
along all the Euler-Lagrange extremals $q(t)$, $t \in [a,b]$, of the problem.
We also note that standard ($\alpha = 1$) Noether's conservation laws
are always a sum of products, as we are assuming in
\eqref{eq:somaPrd}.
\end{rem}

\begin{defn}[\cite{tncdf}]
\label{def:invadf}
Functional \eqref{Pf} is said to be invariant
under the one-parameter group of infinitesimal transformations
\begin{equation}
\label{eq:tinf2}
\begin{cases}
\bar{t} = t + \varepsilon\tau(t,q) + o(\varepsilon) \, ,\\
\bar{q}(\bar{t}) = q(t) + \varepsilon\xi(t,q) + o(\varepsilon) \, ,\\
\end{cases}
\end{equation}
if, and only if,
\begin{equation}
\label{eq:invdf}
\int_{t_{a}}^{t_{b}} L\left(t,q(t),{_{t_a}D_t^\alpha q(t)}\right) dt \\
= \int_{\bar{t}(t_a)}^{\bar{t}(t_b)} L\left(\bar{t},\bar{q}(\bar{t}),
{_{\bar{t}_a}D_{\bar{t}}^\alpha \bar{q}(\bar{t})}\right) d\bar{t}
\end{equation}
for any subinterval $[{t_{a}},{t_{b}}] \subseteq [a,b]$.
\end{defn}

\begin{rem}
Having in mind that condition \eqref{eq:invdf}
is to be satisfied for any subinterval $[{t_{a}},{t_{b}}] \subseteq [a,b]$,
we can rid off the integral signs in \eqref{eq:invdf}
(\textrm{cf.} Definition~\ref{def:inv:gt}).
\end{rem}

The next theorem provides the extension of Noether's theorem
for Fractional Problems of the Calculus of Variations.

\begin{thm}[\cite{tncdf}]
\label{theo:tndf} If functional \eqref{Pf} is invariant
under \eqref{eq:tinf2}, then
\begin{equation*}
\left[L\left(t,q,{_aD_t^\alpha q}\right)
- \alpha\partial_{3} L\left(t,q,{_aD_t^\alpha q}\right)
\cdot{_aD_t^\alpha q} \right] \tau(t,q) \\
+ \partial_{3} L\left(t,q,{_aD_t^\alpha q}\right) \cdot \xi(t,q)
\end{equation*}
is a fractional conservation law (\textrm{cf.} Definition~\ref{eq:fcl}).
\end{thm}


\section{Main Result}

Using Theorem~\ref{theo:tndf}, we obtain here
a formulation of Noether's Theorem for
the fractional optimal control problems
introduced in \cite{Agrawal:2004a}:
\begin{gather}
\label{eq:JO}
I[q(\cdot),u(\cdot)] =\int_a^b L\left(t,q(t),u(t)\right) dt
\longrightarrow \min \, , \\
_aD_t^\alpha q(t)=\varphi\left(t,q(t),u(t)\right) \, , \notag
\end{gather}
together with some boundary conditions on $q(\cdot)$
(which are not relevant with respect to Noether's theorem).
In problem \eqref{eq:JO}, the Lagrangian $L :
[a,b]\times \mathbb{R}^{n} \times \mathbb{R}^{m} \rightarrow
\mathbb{R}$ and the velocity vector $\varphi : [a,b] \times
\mathbb{R}^{n} \times \mathbb{R}^{m} \rightarrow\mathbb{R}^n$ are
assumed to be $C^{1}$ functions with respect to all the arguments.
In agreement with the calculus of variations, we also assume that the
admissible control functions take values on an open set of $\mathbb{R}^m$.

\begin{defn}
A pair $(q(\cdot),u(\cdot))$ satisfying the fractional control
system $_aD_t^\alpha q(t)=\varphi\left(t,q(t),u(t)\right)$ of
problem \eqref{eq:JO}, $t \in [a,b]$, is called a \emph{process}.
\end{defn}

\begin{thm}[\textrm{cf.} (13)-(15) of \cite{Agrawal:2004a}]
\label{th:P}
If $(q(\cdot),u(\cdot))$ is an optimal process for
problem \eqref{eq:JO}, then there exists a co-vector function
$p(\cdot)$ such that the following conditions hold:
\begin{itemize}
\item the Hamiltonian system
\begin{equation*}
\label{eq:Ham}
\begin{cases}
_aD_t^\alpha q(t)&=\partial_4 {\cal H}(t, q(t), u(t),p(t)) \, , \\
_tD_b^\alpha p(t) &= \partial_2{\cal H}(t,q(t),u(t), p(t)) \, ;
\end{cases}
\end{equation*}
\item the stationary condition
\begin{equation*}
\label{eq:CE}
 \partial_3 {\cal H}(t, q(t), u(t), p(t))=0 \, ;
\end{equation*}
\end{itemize}
with the Hamiltonian ${\cal H}$ defined by
\begin{equation}
\label{eq:H} {\cal H}\left(t,q,u,p\right) = L\left(t,q,u\right)
+ p \cdot \varphi\left(t,q,u\right) \, .
\end{equation}
\end{thm}

\begin{rem}
In classical mechanics, the Lagrange multiplier $p$ is called the \emph{generalized momentum}.
In the language of optimal control, $p$ is known as the \emph{adjoint variable}.
\end{rem}

\begin{defn}
\label{def:extPont} Any triplet $(q(\cdot),u(\cdot),p(\cdot))$
satisfying the conditions of Theorem~\ref{th:P} will be called a
\emph{fractional Pontryagin extremal}.
\end{defn}

For the fractional problem of the calculus of variations \eqref{Pf}
one has $\varphi(t,q,u)=u \Rightarrow {\cal H} = L + p \cdot u$,
and we obtain from Theorem~\ref{th:P} that
\begin{gather*}
_aD_t^\alpha q = u \, ,\\
_tD_b^\alpha  p  = \partial_2 L \, ,\\
 \partial_3 {\cal H} = 0 \Leftrightarrow
 p= - \partial_3 L \Rightarrow
{_tD_b}^\alpha p= - _tD_b^\alpha \partial_3 L \, .
\end{gather*}
Comparing the two expressions for $_tD_b^\alpha p$, one arrives to the
Euler-Lagrange differential equations \eqref{eq:eldf}:
$\partial_2 L = - _tD_b^\alpha \partial_3 L$.

We define the notion of invariance for problem \eqref{eq:JO}
in terms of the Hamiltonian, by introducing the augmented functional
as in \cite{Agrawal:2004a}:
\begin{equation}
\label{eq:J} J[q(\cdot),u(\cdot),p(\cdot)] = \\
\int_a^b \left[{\cal
H}\left(t,q(t),u(t),p(t)\right)-p(t) \cdot {_aD_t}^\alpha q(t)\right]dt \, ,
\end{equation}
where ${\cal H}$ is given by \eqref{eq:H}.

\begin{rem}
Theorem~\ref{th:P} is obtained applying the necessary optimality
condition \eqref{eq:eldf} to problem \eqref{eq:J}.
\end{rem}

\begin{defn}
\label{def:inv:gt} A fractional optimal control problem \eqref{eq:JO} is
said to be invariant under the $\varepsilon$-parameter local group
of transformations
\begin{equation}
\label{eq:trf:inf}
\begin{cases}
\bar{t} = t+\varepsilon\tau(t, q(t), u(t), p(t)) + o(\varepsilon) \, , \\
\bar{q}(\bar{t}) = q(t)+\varepsilon\xi(t, q(t), u(t), p(t)) + o(\varepsilon) \, , \\
\bar{u}(\bar{t}) = u(t)+\varepsilon\sigma(t, q(t), u(t), p(t)) + o(\varepsilon) \, , \\
\bar{p}(\bar{t}) = p(t)+\varepsilon \zeta(t, q(t), u(t), p(t))+ o(\varepsilon) \, , \\
\end{cases}
\end{equation}
if, and only if,
\begin{multline}
\label{eq:condInv}
\left[{\cal H}(\bar{t},\bar{q}(\bar{t}),\bar{u}(\bar{t}),\bar{p}(\bar{t}))
-\bar{p}(\bar{t}) \cdot  {_{\bar{a}}D_{\bar{t}}}^\alpha \bar{q}(\bar{t})\right] d\bar{t} \\
=\left[{\cal H}(t,q(t),u(t),p(t))-p(t) \cdot {_aD_t}^\alpha q(t)\right] dt \, .
\end{multline}
\end{defn}

\begin{thm}[Fractional Noether's theorem]
\label{thm:mainResult:FDA06}
If the fractional optimal control problem \eqref{eq:JO} is invariant
under \eqref{eq:trf:inf}, then
\begin{equation}
\label{eq:tndf:CO}
\left[ {\cal H} - \left(1 - \alpha\right) p(t)
\cdot {_aD_t}^\alpha q(t) \right] \tau - p(t) \cdot \xi
\end{equation}
is a fractional conservation law, that is,
\begin{equation*}
\mathcal{D}\left\{
\left[ {\cal H} - \left(1 - \alpha\right) p(t)
\cdot {_aD_t}^\alpha q(t) \right] \tau - p(t) \cdot \xi
\right\} = 0
\end{equation*}
along all the fractional Pontryagin extremals.
\end{thm}

\begin{rem}
For $\alpha = 1$ the fractional optimal control problem \eqref{eq:JO}
is reduced to the classical optimal control problem
\begin{gather*}
I[q(\cdot),u(\cdot)] =\int_a^b L\left(t,q(t),u(t)\right) dt
\longrightarrow \min \, , \\
\dot{q}(t)=\varphi\left(t,q(t),u(t)\right) \, ,
\end{gather*}
and we obtain from Theorem~\ref{thm:mainResult:FDA06}
the optimal control version of Noether's theorem \cite{delfimEJC}:
invariance under a one-parameter group of transformations
\eqref{eq:trf:inf} imply that
\begin{equation}
\label{eq:H9} C(t,q,u,p)={\cal H}(t,q,u,p)\tau-p\cdot \xi
\end{equation}
is constant along any Pontryagin extremal
(one obtains \eqref{eq:H9} from \eqref{eq:tndf:CO}
setting $\alpha = 1$).
\end{rem}

\begin{pf}
The fractional conservation law \eqref{eq:tndf:CO}
is obtained applying Theorem~\ref{theo:tndf}
to the augmented functional \eqref{eq:J}.
\end{pf}


\section{An Example}

Let us consider the autonomous fractional optimal control problem,
\textrm{i.e.} the particular situation when the Lagrangian $L$
and the fractional velocity vector $\varphi$
do not depend explicitly on time $t$:
\begin{equation}
\label{eq:FOCP:CA}
\begin{gathered}
I[q(\cdot),u(\cdot)] =\int_a^b L\left(q(t),u(t)\right) dt
\longrightarrow \min \, , \\
_aD_t^\alpha q(t)=\varphi\left(q(t),u(t)\right) \, .
\end{gathered}
\end{equation}
For the autonomous fractional problem \eqref{eq:FOCP:CA}
the Hamiltonian ${\cal H}$ does not depend explicitly on time,
and it is a simple exercise to check that \eqref{eq:FOCP:CA}
is invariant under time-translations: invariance condition \eqref{eq:condInv}
is satisfied with $\bar{t} = t + \varepsilon$, $\bar{q}(\bar{t}) = q(t)$,
$\bar{u}(\bar{t}) = u(t)$ and $\bar{p}(\bar{t}) = p(t)$. In fact,
given that $d\bar{t} = dt$, \eqref{eq:condInv} holds trivially
proving that ${_{\bar{a}}D_{\bar{t}}}^\alpha \bar{q}(\bar{t})
= {_aD_t}^\alpha q(t)$:
\begin{equation*}
\begin{split}
_{\bar{a}} & D_{\bar{t}}^\alpha \bar{q}(\bar{t}) \\
&= \frac{1}{\Gamma(n-\alpha)}\left(\frac{d}{d\bar{t}}\right)^{n}
\int_{\bar{a}}^{\bar{t}} (\bar{t}-\theta)^{n-\alpha-1}\bar{q}(\theta)d\theta \\
&= \frac{1}{\Gamma(n-\alpha)}\left(\frac{d}{dt}\right)^{n}
\int_{a + \varepsilon}^{t+\varepsilon} (t + \varepsilon-\theta)^{n-\alpha-1}\bar{q}(\theta)d\theta \\
&= \frac{1}{\Gamma(n-\alpha)}\left(\frac{d}{dt}\right)^{n}
\int_{a}^{t} (t-s)^{n-\alpha-1}\bar{q}(s + \varepsilon)ds \\
&= {_{a}D_{t}}^\alpha \bar{q}(t + \varepsilon) = {_{a}D_{t}}^\alpha \bar{q}(\bar{t}) \\
&= {_{a}D_{t}}^\alpha q(t) \, .
\end{split}
\end{equation*}
Using the notation in \eqref{eq:trf:inf} one has
$\tau = 1$ and $\xi = \sigma = \zeta = 0$.
It follows from our fractional Noether's theorem
(Theorem~\ref{thm:mainResult:FDA06}) that
\begin{equation}
\label{eq:ConsHam:alpha}
\mathcal{D}\left\{
{\cal H} - \left(1 - \alpha\right) p(t)
\cdot {_aD_t}^\alpha q(t)
\right\} = 0 \, .
\end{equation}
In the classical framework of optimal control theory $\alpha = 1$
and our operator $\mathcal{D}$ coincides with $\frac{d}{dt}$:
we then get from \eqref{eq:ConsHam:alpha}
the well known fact that the Hamiltonian is a preserved
quantity along any Pontryagin extremal.


\section{Conclusions}
\label{sec:Conc}

The proof of fractional Euler-Lagrange equations
is a subject of strong current study
\cite{CD:Riewe:1996,CD:Riewe:1997,Klimek2001,CD:Agrawal:2002,%
MR1966935,CD:BalAv:2004,Klimek2005,El-Nabulsi2005b,El-Nabulsi2005a}
because of its numerous applications.
In \cite{tncdf} a fractional Noether's theorem is proved.

The fractional variational theory is in its childhood
so that much remains to be done. This is particularly true
in the area of fractional optimal control,
where the results are rare. A fractional
Hamiltonian formulation is obtained in \cite{Muslih:2005},
but only for systems with linear velocities.
The main study of fractional optimal control problems
seems to be \cite{Agrawal:2004a}, where the Euler-Lagrange
equations for fractional optimal control problems
(Theorem~\ref{th:P}) are obtained,
using the traditional approach of the Lagrange multiplier rule.
Here we use the Lagrange multiplier technique to derive,
from the results in \cite{tncdf}, a Noether-type theorem
for fractional optimal control systems. As an example
we have considered a fractional autonomous problem,
proving that the Hamiltonian defines a conservation law only
in the case $\alpha = 1$.


\section*{Acknowledgment}

We are grateful to Professor Tenreiro Machado
for drawing our attention to the \emph{2nd IFAC Workshop on Fractional
Differentiation and its Applications},
19--21 July, 2006, Porto, Portugal,
and for encouraging us to write the present work.



\end{document}